\UseRawInputEncoding
\documentclass[11pt]{article}
\usepackage{amsfonts}
\usepackage{color}
\usepackage{amsmath,amssymb,comment}
\usepackage{verbatim}
\newtheorem{theo}{Theorem}[section]
\newtheorem{lem}[theo]{Lemma}

\newcommand{\mysection}[1]{\section{#1} \setcounter{equation}{0}}
    \makeatletter
\def\@fnsymbol#1{\ensuremath{\ifcase#1\or *\or \ddagger\or
   \mathsection\or \mathparagraph\or \|\or **\or \dagger\dagger
   \or \ddagger\ddagger \else\@ctrerr\fi}}
    \makeatother
\newcommand{\proof}{{\sc Proof.} \quad}
\newcommand{\proofc}{{\sc Proof} \ }
\newcommand{\be}{\begin{equation} \label}
\newcommand{\ee}{\end{equation}}
\newcommand{\bea}{\begin{eqnarray}\label}
\newcommand{\eea}{\end{eqnarray}}
\newcommand{\bas}{\begin{eqnarray*}}
\newcommand{\eas}{\end{eqnarray*}}
\newcommand{\bit}{\begin{itemize}}
\newcommand{\eit}{\end{itemize}}
\newcommand{\qed}{\hfill$\Box$ \vskip.2cm}
\newcommand{\nn}{\nonumber}
\newcommand{\R}{\mathbb{R}}
\newcommand{\N}{\mathbb{N}}
\newcommand{\pO}{\partial\Omega}

\newcommand{\eps}{\varepsilon}

\newcommand{\wto}{\rightharpoonup}
\newcommand{\wsto}{\stackrel{\star}{\rightharpoonup}}
\newcommand{\hra}{\hookrightarrow}
\newcommand{\io}{\int_\Omega}
\newcommand{\na}{\nabla}
\newcommand{\Del}{\Delta}
\newcommand{\del}{\delta}

\newcommand{\lam}{\lambda}
\newcommand{\sig}{\sigma}
\newcommand{\pa}{\partial}
\newcommand{\bom}{\overline{\Omega}}
\newcommand{\Om}{\Omega}

\newcommand{\ov}{\overline}

\newcommand{\hs}{\hspace*}
\newcommand{\sm}{\setminus}

\newcommand{\vp}{\varphi}
\newcommand{\lbal}{\left\{ \begin{array}{l}}
\newcommand{\lball}{\left\{ \begin{array}{ll}}
\newcommand{\ear}{\end{array} \right.}
\newcommand{\ouz}{\ov{u}_0}

\newcommand{\abs}{\\[5pt]}

\newcommand{\oy}{\ov{y}}

\newcommand{\ueps}{u_\eps}
\newcommand{\veps}{v_\eps}
\newcommand{\weps}{w_\eps}
\newcommand{\heps}{h_\eps}

\newcommand{\geps}{g_\eps}
\newcommand{\yeps}{y_\eps}

\begin{document}
\enlargethispage{10mm}
\title{Refined regularity analysis for a Keller-Segel-consumption system\\
involving signal-dependent motilities}
\author{
Genglin Li\footnote{1182028@mail.dhu.edu.cn}\\
{\small College of Information Science and Technology, }\\
{\small Donghua University, Shanghai 201620, P.R.~China}
\and
Michael Winkler\footnote{michael.winkler@math.uni-paderborn.de}\\
{\small Institut f\"ur Mathematik, Universit\"at Paderborn,}\\
{\small 33098 Paderborn, Germany} 
}
\date{}
\maketitle
\begin{abstract}
\noindent 
We consider the Keller-Segel-type migration-consumption system involving signal-dependent motilities,
\bas
	\lbal
	u_t = \Delta \big(u\phi(v)\big), \\[1mm]
	v_t = \Delta v-uv,
	\ear
\eas
in smoothly bounded domains $\Om\subset\R^{n}$, $n\ge 1$.
Under the assumption that $\phi\in C^1([0,\infty))$ is positive on $[0,\infty)$, and for nonnegative initial data from 
$(C^0(\bom))^\star \times L^\infty(\Om)$,
previous literature has provided results on global existence of certain very weak solutions with possibly quite poor
regularity properties, and on large time stabilization toward semitrivial equilibria with respect to the topology in 
$(W^{1,2}(\Om))^\star \times L^\infty(\Om)$.\abs
The present study reveals that solutions in fact enjoy significantly stronger regularity features when
$0<\phi\in C^3([0,\infty))$ and the initial data belong to $(W^{1,\infty}(\Om))^2$:
It is firstly shown, namely, that then in the case $n\le 2$
an associated no-flux initial-boundary value problem even admits a global classical solution, and that each of these solutions smoothly stabilizes in the sense that as $t\to\infty$ we have 
\bas
	u(\cdot,t) \to \frac{1}{|\Om|} \io u_0
	\qquad \mbox{and} \qquad
	v(\cdot,t)\to 0
	\qquad \qquad (\star)
\eas
even with respect to the norm in $L^\infty(\Om)$ in both components.\abs
In the case when $n\ge 3$, secondly, some genuine weak solutions are found to exist globally, inter alia satisfying
$\na u\in L^\frac{4}{3}_{loc}(\bom\times [0,\infty);\R^n)$. 
In the particular three-dimensional setting, any such solution is seen to become eventually smooth and to satisfy ($\star$).\abs
\noindent {\bf Key words:} chemotaxis; global existence; boundedness; eventual regularity; large time behavior\\
 {\bf MSC 2020:} 35B65 (primary); 35B40, 35K55, 35Q92, 92C17 (secondary)
\end{abstract}
\newpage
\section{Introduction}\label{intro}
In the macroscopic modeling of bacterial migration influenced by local sensing mechanisms,
a prominent role is played by parabolic systems of the form
\be{01}
 	\left\{ \begin{array}{l}	
 	u_t = \Del \big(u\phi(v)\big), \\[1mm]
 	v_t = \Delta v + f(u,v)
 	\end{array} \right.
\ee
(\cite{fu}, \cite{liu}, \cite{KS3}).
By featuring a particular link between diffusion and cross-diffusion in the migration operator determining its first equation,
(\ref{01}) can be considered a special case of more general chemotaxis systems of Keller-Segel type (\cite{KS1});
indeed, suitably making use of this characteristic structure has been underlying substantial parts of the existing 
literature on subclasses of (\ref{01}).
In the case when $f(u,v)= u-v$, for instance,
within large classes of nonlinearities $\phi$ an
appropriate exploitation of this core property has been forming the respective key in a considerable collection of
studies concerned with basic solution theories for the corresponding chemotaxis-production version of (\ref{01})
(\cite{fujie_senba}, \cite{ahn_yoon}, \cite{fujie_jiang_ACAP2021}, \cite{taowin_M3AS}, \cite{burger}, \cite{desvillettes}),
as well as some close variants
(\cite{jiang_laurencot}, \cite{fujie_jiang_JDE2020}, \cite{jiang_arxiv},
\cite{jin_kim_wang}, \cite{wang_wang}, \cite{yifu_wang}, \cite{wenbin_lv_ZAMP}, \cite{wenbin_lv_EECT});
beyond this, results on spatially structured large time behavior, and even on the occurrence of infinite-time blow-up
(\cite{fujie_jiang_CVPDE}, \cite{fujie_senba}, \cite{jin_wang}, \cite{ahn_yoon}, \cite{DLTW}),
provide conclusive evidence for significant structure-supporting features of such representatives of (\ref{01})
that were already observed in numerical simulations (\cite{desvillettes}).\abs
In comparison to this, the corresponding chemotaxis-consumption counterpart, as given by
\be{00}
 	\left\{ \begin{array}{l}	
 	u_t = \Del \big(u\phi(v)\big), \\[1mm]
 	v_t = \Delta v -uv,
 	\end{array} \right.
\ee
appears to exhibit a significantly weaker tendency toward pattern generation:
Here the underlying hypothesis that the considered signal is consumed by cells, rather than produced, does not only 
go along with an evident nonexistence of inhomogeneous steady states; indeed, as has recently been shown in 
\cite{liwin2}, 
indeed, when determined by strictly positive migration rate functions $\phi\in C^1([0,\infty))$,
the taxis-absorption interaction in (\ref{00}) even enforces certain global very weak solutions,
known to exist for all initial data from $(C^0(\bom))^\star \times L^\infty(\Om)$ (\cite{liwin2}; cf.~also \cite{li_zhao_ZAMP}), 
to approach the semi-trivial
equilibrium $(\frac{1}{|\Om|} \io u_0, 0)$ in the topology of $(W^{1,2}(\Om))^\star \times L^\infty(\Om)$ in the large time
limit when $n\le 3$ (see \cite{liwin2} and also Lemma \ref{lem1} below);
nontrivial large time dynamics in (\ref{00}) seems possible only in degenerate cases in which $\phi(0)=0$
(\cite{win_sig_dep_mot_cons_2}, \cite{win_sig_dep_mot_cons_largetime}).\abs
{\bf Main results.} \quad
The objective of this manuscript is to undertake a refined regularity analysis for solutions of (\ref{00}) in said
non-degenerate setting involving strictly positive $\phi$;
in fact, according to the mild assumptions on regularity of $\phi$ and the initial data made in \cite{liwin2}, 
only quite poor regularity information on the global solutions constructed there seems available so far:
The corresponding first components, for instance, are merely known to belong to 
$L^\infty((0,\infty);L^1(\Om)) \cap L^2_{loc}(\bom\times (0,\infty))$, and to be continuous on $[0,\infty)$ 
as $C^0(\bom))^\star$-valued functions with respect to the weak topology therein;
in particular, this does not rule out the emergence of transient singularities, and due to the absence of
integrability information on $\nabla u$, this regularity class is far too large to
let these objects play the role of genuine weak solutions.\abs
Our results will inter alia reveal that under moderately stronger assumptions than those from \cite{liwin2}, particularly
still not relying on any requirement on $\phi$ that goes beyond smoothness and positivity,
the problem (\ref{00}) actually admits global classical solutions when $n\le 2$, and eventually smooth weak solutions
when $n=3$, and that in both these cases the obtained solutions smoothly stabilize.
In this regard, (\ref{00}) will thus be seen to share essential
features of blow-up exclusion with the related classical chemotaxis-consumption
system with its first equation, $u_t=\Del u - \na\cdot (u\na v)$, reflecting interplay between linear diffusion
and standard Keller-Segel type cross-diffusion (\cite{taowin_consumption}).\abs
Indeed, for the initial-boundary value problem
\be{0}
	\left\{ \begin{array}{ll}	
	u_t = \Del \big(u\phi(v)\big),
	\qquad & x\in\Om, \ t>0, \\[1mm]
	v_t = \Delta v-uv,
	\qquad & x\in\Om, \ t>0, \\[1mm]
	\na \big(u\phi(v)\big) \cdot \nu = \na v \cdot \nu=0,
	\qquad & x\in\pO, \ t>0, \\[1mm]
	u(x,0)=u_0(x), \quad v(x,0)=v_0(x),
	\qquad & x\in\Om,
	\end{array} \right.
\ee
the first of our main results addresses the corresponding low-dimensional framework:
\begin{theo}\label{theo9}
  Let $n\in\{1,2\}$ and $\Om\subset\R^n$ be a bounded domain with smooth boundary, and suppose that
  \be{phi}
	\phi\in C^3([0,\infty))
	\mbox{\quad is such that \quad } \phi(\xi)>0 \mbox{ for all } \xi\ge 0,
  \ee
  and that 		
  \be{init}
	\lbal 
	u_0\in W^{1,\infty}(\Om)
	\mbox{ is nonnegative with $u_0\not\equiv 0$, and} \\[1mm]
	v_0\in W^{1,\infty}(\Om)
	\mbox{ is nonnegative.}
	\ear
  \ee
  Then one can find nonnegative functions
  \be{9.01}
	\lbal
	u\in C^0(\bom\times [0,\infty)) \cap C^{2,1}(\bom\times (0,\infty))
	\qquad \mbox{and} \\[1mm]
	v\in C^0(\bom\times [0,\infty)) \cap C^{2,1}(\bom\times (0,\infty))
	\ear
  \ee
  such that $(u,v)$ solves (\ref{0}) in the classical sense. Moreover,
  \be{9.1}
	u(\cdot,t)\to \frac{1}{|\Om|} \io u_0
	\quad \mbox{in $L^\infty(\Om)$ \quad and \quad}
	v(\cdot,t)\to 0
	\quad \mbox{in $L^\infty(\Om)$ \qquad as } t\to\infty.
  \ee
\end{theo}
Under the same hypotheses, our result concerning the higher-dimensional version of (\ref{0}) can be stated as follows.
\begin{theo}\label{theo17}
  Let $n\ge 3$ and $\Om\subset\R^n$ be a smoothly bounded domain, and assume (\ref{phi}) and (\ref{init}).
  Then one can find
  \be{reg2}
	\lbal
	u\in L^2_{loc}(\bom\times [0,\infty)) \cap L^\frac{4}{3}_{loc}([0,\infty);W^{1,\frac{4}{3}}(\Om))
	\qquad \mbox{and} \\[1mm]
	v\in L^\infty(\Om\times (0,\infty)) \cap L^2_{loc}([0,\infty);W^{2,2}(\Om))
	\ear
  \ee
  such that $u\ge 0$ and $v\ge 0$ a.e.~in $\Om\times (0,\infty)$, and that $(u,v)$ forms a global weak solution of (\ref{0})
  in the sense that
  \be{wu}
	- \int_0^\infty \io u\vp_t - \io u_0\vp(\cdot,0)
	= - \int_0^\infty \io \na\big(u\phi(v)\big)\cdot\na\vp
  \ee
  and
  \be{wv}
	\int_0^\infty \io v\vp_t + \io v_9\vp(\cdot,0)
	= \int_0^\infty \io \na v\cdot\na\vp + \io uv\vp
  \ee
  for all $\vp\in C_0^\infty(\bom\times [0,\infty))$.\abs
  If moreover $n=3$, then there exists $T>0$ such that after re-definition of $(u,v)$ on a null set we can achieve that
  \be{17.1}
	u\in C^{2,1}(\bom\times (T,\infty))
	\qquad \mbox{and} \qquad
	v\in C^{2,1}(\bom\times (T,\infty)),
  \ee
  and that (\ref{9.1}) holds.
\end{theo}
{\bf Main ideas.} \quad
The starting point of our analysis consists in making appropriate use of the non-degenerate migration operator 
contained in the first equation in (\ref{0}) to derive a spatio-temporal a priori estimate for $u$, which is 
accomplished by a duality-based argument (Lemma \ref{lem2}), thereby facilitating the establishment of certain regularity 
information on $v$ (Lemma \ref{lem3}). \abs
In light of these estimates, in the lower-dimensional setting of Theorem \ref{theo9}, 
boundedness of $u$ in $L^p(\Omega)$ for any
$p\ge1$ can be established through a straightforward testing procedure (Lemma \ref{lem4}), providing the first step
in a bootstrap-type argument that involves standard regularity theory for parabolic equations and yields
smoothness of solutions and the stabilization result recorded in Theorem \ref{theo9} (Lemmata \ref{lem5}-\ref{lem8}). \abs
For spatially higher-dimensional versions of (\ref{0}), another application of the estimates obtained 
in Section 2 leads to time-independent bounds for expressions of the form $\io u\ln u$, and to an estimate for $\na u$
(Lemma \ref{lem10}), hence implying
that the limit function $(u,v)$ constructed in \cite{liwin2} actually is a global weak solution as claimed in Theorem \ref{theo17}.
Furthermore, in the three-dimensional case
an eventual smallness property of $\|v\|_{L^\infty(\Om)}$ (Lemma \ref{lem12})
facilitates an analysis of functionals of the form
$\io \frac{u^p(\cdot,t)}{(\delta-v(\cdot,t))^{\kappa}}$, with arbitrary $p>1$ and suitably chosen $\del=\del(p)>0$ 
and $\kappa=\kappa(p)>0$, to establish $L^p$ bounds for $u$ (Lemma \ref{lem13}).
Subsequent higher order regularity analysis
(Lemmata \ref{lem14}-\ref{lem16}) will lead to the conclusion of Theorem \ref{theo17}.
\mysection{Preliminaries}
In order to appropriately regularize (\ref{0}), let us consider the approximate variants of (\ref{0}) given by 
\be{0eps}
	\lball
	u_{\eps t} = \Del \big(\ueps\phi(\veps)\big),
	\qquad & x\in\Omega, \ t>0, \\[1mm]
	v_{\eps t} = \Delta\veps - \frac{\ueps\veps}{1+\eps\ueps},
	\qquad & x\in\Omega, \ t>0, \\[1mm]
	\frac{\partial\ueps}{\partial\nu}=\frac{\partial\veps}{\partial\nu}=0,
	\qquad & x\in\pO, \ t>0, \\[1mm]
	\ueps(x,0)=u_0(x), \quad \veps(x,0)=v_0(x),
	\qquad & x\in\Om,
	\ear
\ee
for $\eps\in (0,1)$. 
Then in the context of actually less restrictive assumptions on $\phi$ and the initial data, the following
collection of findings extracts from \cite[Theorems 1.1 and 1.2 and Lemma 2.2]{liwin2} 
what will be needed for our subsequent analysis.
In formulating this, given a Banach space $X$ we let $C^0_{w-\star}([0,\infty); X)$ 
denote the space of functions that are continuous on $[0,\infty)$ with respect to the weak-$\star$ topology in $X$.
\begin{lem}\label{lem1}
  Let $\eps\in (0,1)$. Then there exist nonnegative functions
  \be{01.1}
	\lbal
	\ueps\in C^0(\bom\times [0,\infty)) \cap C^{2,1}(\bom\times (0,\infty))
	\qquad \mbox{and} \\[1mm]
	\veps \in \bigcap_{q>n} C^0([0,\infty);W^{1,q}(\Om)) \cap C^{2,1}(\bom\times (0,\infty))
	\ear
  \ee
  such that $(\ueps,\veps)$ solves (\ref{0eps}) in the classical sense, and that
  \be{mass}
  	\io \ueps(\cdot,t) = \io u_0
  	\qquad \mbox{for all $t>0$}
  \ee
  as well as
  \be{vinfty}
  	\|\veps(\cdot,t)\|_{L^\infty(\Om)} \le \|\veps(\cdot,t_0)\|_{L^\infty(\Om)}
  	\qquad \mbox{for all $t_0\ge 0$ and $t>t_0$.}
  \ee
  Moreover, there exist $(\eps_j)_{j\in\N}\subset (0,1)$ as well as nonnegative functions
  \be{regw}
	\lbal
	u\in C^0_{w-\star}([0,\infty);(C^0(\bom))^\star) \cap L^\infty((0,\infty);L^1(\Om)) \cap L^2_{loc}(\bom\times (0,\infty))
	\qquad \mbox{and} \\[1mm]
	v\in C^0_{w-\star}([0,\infty;L^\infty(\Om)) \cap L^\infty(\Om\times (0,\infty)) \\
	\hs{10mm}
	\cap L^2_{loc}((0,\infty);W^{2,2}(\Om)) \cap L^4_{loc}((0,\infty);W^{1,4}(\Om)) \cap L^\infty_{loc}((0,\infty);W^{1,2}(\Om))
	\ear
  \ee
  such that $\eps_j\searrow 0$ as $j\to\infty$, that
  \be{1.1}
	\ueps\to u
	\quad \mbox{and} \quad
	\veps\to v
	\qquad \mbox{a.e.~in } \Om\times (0,\infty)
  \ee
  as $\eps=\eps_j\searrow 0$, and that $(u,v)$ forms a global very weak solution of (\ref{0}) in the sense that
  $u(\cdot,0)=u_0$ in $(C^0(\bom))^\star$ and $v(\cdot,0)=v_0$ in $L^\infty(\Om)$,
  and that 
  for each $\vp\in C_0^\infty(\bom\times (0,\infty))$ fulfilling $\frac{\pa\vp}{\pa\nu}=0$ on $\pO\times (0,\infty)$ we have
  \bea{wu1}
	-\int_0^\infty \io u\vp_t
	= \int_0^\infty \io u\phi(v) \Del\vp
  \eea
  and
  \be{wv1}
	- \int_0^\infty \io v\vp_t  
	= \int_0^\infty \io v\Del\vp
	- \int_0^\infty \io uv\vp.
  \ee
  Finally, if $n\le 3$ then there exists a null set $N\subset (0,\infty)$ such that
  \be{1.2}
	u(\cdot,t)\to \frac{1}{|\Om|} \io u_0
	\ \mbox{in } (W^{1,2}(\Om))^\star
	\quad \mbox{and} \quad
	v(\cdot,t)\to 0
	\ \mbox{in } L^\infty(\Om)
	\qquad \mbox{as } (0,\infty)\sm N \ni t\to\infty.
  \ee
\end{lem}
A yet fairly elementary but crucial regularity feature of the first solution component can be obtained by means of 
a duality-based reasoning in the spirit of a precedent contained, e.g., in \cite{taowin_M3AS}.
\begin{lem}\label{lem2}
  There exists $C>0$ such that
  \be{2.1}
	\int_{(t-1)_+}^t \io \ueps^2 \le C
	\qquad \mbox{for all $t>0$ and } \eps\in (0,1).
  \ee
\end{lem}
\proof
  Following \cite{taowin_M3AS}, we let $A$ denote the realization of $-\Del+1$ under homogeneous Neumann boundary conditions
  in $L^2(\Om)$, and multiply the identity
  $\pa_t A^{-1} \ueps + \ueps\phi(\veps)=A^{-1}(\ueps\phi(\veps))$ by $\ueps$ to see that by self-adjointness of $A^{-1}$,
  and by Young's inequality,
  \bas
	\frac{1}{2} \frac{d}{dt} \io \big|A^{-\frac{1}{2}}\ueps\big|^2
	+ \io \ueps^2 \phi(\veps)
	&=& \io \ueps A^{-1}\big(\ueps\phi(\veps)\big) \\
	&=& \io \ueps\phi(\veps) A^{-1} \ueps \\
	&\le& \frac{1}{2} \io \ueps^2 \phi(\veps)
	+ \frac{1}{2} \io \phi(\veps) \big|A^{-1} \ueps\big|^2
	\qquad \mbox{for all $t>0$ and } \eps\in (0,1).
  \eas
  Since (\ref{phi}) together with (\ref{vinfty}) yields $c_1>0$ and $c_2>0$ such that $c_1\le \phi(\veps) \le c_2$ in 
  $\Om\times (0,\infty)$ for all $\eps\in (0,1)$, this implies that
  $\yeps(t):=\io |A^{-\frac{1}{2}}\ueps(\cdot,t)|^2$, $t\ge 0, \eps\in (0,1)$, satisfies
  \be{2.2}
	\yeps'(t) + \yeps(t)
	+ c_1 \io \ueps^2 
	\le c_2 \io \big|A^{-1}\ueps\big|^2
	+ \io \big|A^{-\frac{1}{2}}\ueps\big|^2
	\qquad \mbox{for all $t>0$ and } \eps\in (0,1).
  \ee
  Here since the boundedness of $\Om$ ensures that $L^2(\Om)$ is compactly embedded into both domains $D(A^{-1})$ 
  and $D(A^{-\frac{1}{2}})$, by means of associated Ehrling inequalities and (\ref{mass}) we readily find $c_3>0$ and $c_4>0$
  such that
  \bas
	c_2 \io \big| A^{-1}\ueps\big|^2 
	+ \io \big|A^{-\frac{1}{2}}\ueps\big|^2
	\le \frac{c_1}{2} \io \ueps^2
	+ c_3 \cdot \bigg\{ \io \ueps\bigg\}^2
	\le \frac{c_1}{2} \io \ueps^2 + c_4
	\quad \mbox{for all $t>0$ and } \eps\in (0,1),
  \eas
  so that from (\ref{2.2}) we obtain that
  \bas
	\yeps'(t) + \yeps(t) + \frac{c_1}{2} \io \ueps^2 \le c_4
	\qquad \mbox{for all $t>0$ and } \eps\in (0,1).
  \eas
  Therefore, $\yeps(t)\le c_5:=\max\big\{ \io |A^{-\frac{1}{2}} u_0|^2 \, , \, c_4\big\}$ for all $t\ge 0$ and $\eps\in (0,1)$, and 
  thus
  \bas
	\frac{c_1}{2} \int_{(t-1)_+}^t \io \ueps^2
	\le \yeps \big( (t-1)_+\big) + c_4 \le c_5+c_4
  \eas
  for all $t>0$ and $\eps\in (0,1)$.
\qed
As a first application, when combined with a standard testing procedure the above yields some basic regularity information
on $\veps$, as the previous estimate available without any restriction on the spatial dimension.
\begin{lem}\label{lem3}
  There exists $C>0$ such that
  \be{3.01}
	\io \big|\na\veps(\cdot,t)\big|^2 \le C
	\qquad \mbox{for all $t>0$ and } \eps\in (0,1)
  \ee
  and
  \be{3.1}
	\int_{(t-1)_+}^t \io |\Del\veps|^2 \le C
	\qquad \mbox{for all $t>0$ and } \eps\in (0,1)
  \ee
  as well as
  \be{3.2}
	\int_{(t-1)_+}^t \io |\na\veps|^4 \le C
	\qquad \mbox{for all $t>0$ and } \eps\in (0,1)
  \ee
  and
  \be{3.22}
	\int_{(t-1)_+}^t \io v_{\eps t}^2 \le C
	\qquad \mbox{for all $t>0$ and } \eps\in (0,1).
  \ee
\end{lem}
\proof
  According to the second equation in (\ref{0eps}), Young's inequality and (\ref{vinfty}), 
  \bas
	\frac{1}{2} \frac{d}{dt} \io |\na\veps|^2
	+ \io |\Del\veps|^2
	&=& \io \frac{\ueps\veps}{1+\eps\ueps} \Del\veps \\
	&\le& \frac{1}{2} \io |\Del\veps|^2
	+ \frac{1}{2} \io \Big(\frac{\ueps\veps}{1+\eps\ueps}\Big)^2 \\
	&\le& \frac{1}{2} \io |\Del\veps|^2
	+ \frac{c_1^2}{2} \io \ueps^2
	\qquad \mbox{for all $t>0$ and } \eps\in (0,1),
  \eas
  where $c_1:=\|v_0\|_{L^\infty(\Om)}$.
  Since the Gagliardo-Nirenberg inequality combined with elliptic regularity theory provides $c_2>0$ such that again by 
  (\ref{vinfty}) we have
  \bas
	\io |\na\veps|^4 \le c_2\|\Del\veps\|_{L^2(\Om)}^2 \|\veps\|_{L^\infty(\Om)}^2
	\le c_1^2 c_2 \io |\Del\veps|^2
	\qquad \mbox{for all $t>0$ and } \eps\in (0,1),
  \eas
  and since thus also
  \bas
	\io |\na\veps|^2 
	\le \frac{1}{4} \io |\Del\veps|^2
	+ c_1^2 c_2 |\Om|
	\qquad \mbox{for all $t>0$ and } \eps\in (0,1)
  \eas
  thanks to Young's inequality, this implies that
  \bas
	\yeps(t):=\io \big|\na\veps(\cdot,t)\big|^2,
	\qquad t\ge 0, \eps\in (0,1),
  \eas
  as well as
  \bas
	\geps(t):=\frac{1}{2} \io \big|\Del\veps(\cdot,t)\big|^2
	+ \frac{1}{4c_1^2 c_2} \io \big|\na\veps(\cdot,t)\big|^4
	\ \mbox{and} \
	\heps(t):=c_1^2 \io \ueps^2(\cdot,t)
	+ c_1^2 c_2 |\Om|,
	\quad t\ge 0, \eps\in (0,1),
  \eas
  satisfy
  \bea{3.3}
	\yeps'(t) + \yeps(t)
	&\le& \bigg\{ - \io |\Del\veps|^2 + c_1^2 \io \ueps^2\bigg\}	
	+ \bigg\{ \frac{1}{4} \io |\Del\veps|^2
	+ c_1^2 c_2 |\Om| \bigg\} \nn\\
	&\le& - \frac{1}{2} \io |\Del\veps|^2
	- \frac{1}{4} \io |\Del\veps|^2
	+ \heps(t) \nn\\
	&\le& - \geps(t) + \heps(t)
	\qquad \mbox{for all $t>0$ and } \eps\in (0,1).
  \eea
  As from Lemma \ref{lem2} we obtain $c_3>0$ such that
  \be{3.4}
	\int_{(t-1)_+}^t \heps(s) ds \le c_3
	\qquad \mbox{for all $t>0$ and } \eps\in (0,1),
  \ee
  using that
  \bas
	\int_0^t e^{-\lam (t-s)} \heps(s) ds
	\le \frac{1}{1-e^{-\lam}} \cdot \sup_{s>0} \int_{(s-1)_+}^s \heps(\sig)d\sig
	\qquad \mbox{for all $t>0$, $\eps\in (0,1)$ and } \lam>0
  \eas
  (cf.~\cite[Lemma 3.4]{win_JFA}), from (\ref{3.3}) we infer that $\yeps(t)\le c_4:=\io \big|\na v_0\big|^2 + \frac{c_3}{1-e^{-1}}$
  for all $t\ge 0$ and $\eps\in (0,1)$, and that hence, again by (\ref{3.4}),
  \bas
	\int_{(t-1)_+}^t \geps(s) ds
	\le \yeps \big( (t-1)_+\big) + \int_{(t-1)_+}^t \heps(s) ds
	\le c_4 + c_3
	\qquad \mbox{for all $t>0$ and } \eps\in (0,1).
  \eas
  In view of the definitions of $(\yeps)_{\eps\in (0,1)}$ and $(\geps)_{\eps\in (0,1)}$, these inequalities establish both
  (\ref{3.01}) and (\ref{3.1})-(\ref{3.2}), whereupon (\ref{3.22}) immediately results upon estimating
  $v_{\eps t}^2 \le 2|\Del\veps|^2 + 2\ueps^2\veps^2$ for $\eps\in (0,1)$, and combining (\ref{3.1}) with
  Lemma \ref{lem2} and (\ref{vinfty}).
\qed
\mysection{The case $n\le 2$. Proof of Theorem \ref{theo9}}
In this section the global classical solvability for problem (\ref{0}) as well as large time convergence of solutions is examined 
in the case when $n\le 2$. 
The following lemma already contains the main step toward this, establishing $L^p$ bounds for $\ueps$ on the basis
of the spatio-temporal $L^4$ estimate for $\na\veps$ provided by (\ref{3.2}):
\begin{lem}\label{lem4}
  Let $n\le 2$. Then for all $p>1$ there exists $C(p)>0$ such that
  \be{4.1}
	\io \ueps^p(\cdot,t) \le C(p)
	\qquad \mbox{for all $t>0$ and } \eps\in (0,1).
  \ee
\end{lem}
\proof
  We detail the proof only for $n=2$, and remark that the one-dimensional case can be addressed by very minor modification
  of the argument.\\
  Using the first equation in (\ref{0eps}) along with Young's inequality, we then see that
  \bas
	\frac{1}{p} \frac{d}{dt} \io \ueps^p
	&=& -(p-1) \io \ueps^{p-2} \phi(\veps) |\na\ueps|^2
	- (p-1) \io \ueps^{p-1} \phi'(\veps) \na\ueps\cdot\na\veps \\
	&\le& -\frac{p-1}{2} \io \ueps^{p-2} \phi(\veps) |\na\ueps|^2
	+ \frac{p-1}{2} \io \ueps^p \frac{\phi'^2(\veps)}{\phi(\veps)} |\na\veps|^2
	\qquad \mbox{for all $t>0$ and } \eps\in (0,1),
  \eas
  whence again combining (\ref{phi}) with (\ref{vinfty}) we obtain $c_1=c_1(p)>0$ and $c_2=c_2(p)>0$ such that
  \be{4.2}
	\frac{d}{dt} \io \ueps^p
	+ c_1 \io \big|\na\ueps^\frac{p}{2}\big|^2
	\le c_2 \io \ueps^p |\na\veps|^2
	\qquad \mbox{for all $t>0$ and } \eps\in (0,1).
  \ee
  Here 		
  by the Cauchy-Schwarz inequality, the Gagliardo-Nirenberg inequality, (\ref{mass}) and Young's inequality we find
  that with some positive constants $c_i=c_i(p), i\in\{3,4,5,6\}$, we have
  \bea{4.3}
	c_2 \io \ueps^p |\na\veps|^2
	&\le& c_2\|\na\veps\|_{L^4(\Om)}^2 \|\ueps^\frac{p}{2}\|_{L^4(\Om)}^2 \nn\\
	&\le& c_3 \|\na\veps\|_{L^4(\Om)}^2 \|\na\ueps^\frac{p}{2}\|_{L^2(\Om)} \|\ueps^\frac{p}{2}\|_{L^2(\Om)}
	+ c_3 \|\na\veps\|_{L^4(\Om)}^2 \|\ueps^\frac{p}{2}\|_{L^\frac{2}{p}(\Om)}^2 \nn\\
	&\le& c_3 \|\na\veps\|_{L^4(\Om)}^2 \|\na\ueps^\frac{p}{2}\|_{L^2(\Om)} \|\ueps^\frac{p}{2}\|_{L^2(\Om)}
	+ c_4 \|\na\veps\|_{L^4(\Om)}^2 \nn\\
	&\le& \frac{c_1}{2} \io \big|\na\ueps^\frac{p}{2}\big|^2
	+ c_5 \|\na\veps\|_{L^4(\Om)}^4 \io \ueps^p
	+ c_4 \|\na\veps\|_{L^4(\Om)}^2 \nn\\
	&\le& \frac{c_1}{2} \io \big|\na\ueps^\frac{p}{2}\big|^2
	+ c_6 \cdot \bigg\{ \io |\na\veps|^4 +1\bigg\} \cdot \bigg\{ \io \ueps^p+1\bigg\}
  \eea
  for all $t>0$ and $\eps\in (0,1)$.
  Since we can similarly apply Young's inequality, the Gagliardo-Nirenberg inequality and (\ref{mass}) to fix
  $c_7=c_7(p)>0$ and $c_8=c_8(p)>0$ such that
  \bas
	\bigg\{ \io \ueps^p +1 \bigg\}^\frac{p}{p-1}
	&\le& 2^\frac{1}{p-1} \|\ueps^\frac{p}{2}\|_{L^2(\Om)}^\frac{2p}{p-1} + 2^\frac{1}{p-1} \\
	&\le& c_7 \|\na\ueps^\frac{p}{2}\|_{L^2(\Om)}^2 \|\ueps^\frac{p}{2}\|_{L^\frac{2}{p}(\Om)}^\frac{2}{p-1}
	+ c_7 \|\ueps^\frac{p}{2}\|_{L^\frac{2}{p}(\Om)}^\frac{2p}{p-1} + 2^\frac{1}{p-1} \\
	&\le& c_8 \io \big|\na\ueps^\frac{p}{2}\big|^2 + c_8
	\qquad \mbox{for all $t>0$ and } \eps\in (0,1),
  \eas
  from (\ref{4.2}) and (\ref{4.3}) we thus infer that for
  \bas
	\yeps(t):=\io \ueps^p(\cdot,t)+1,
	\ t\ge 0, \eps\in (0,1),
	\quad \mbox{and} \quad
	\heps(t):=c_6 \io \big|\na\veps(\cdot,t)\big|^4 + c_6 + \frac{c_1}{2},
	\ t>0, \eps\in (0,1),
  \eas
  we have
  \bas
	\yeps'(t) + \frac{c_1}{2c_8} \yeps^\frac{p}{p-1}(t) \le \heps(t) \yeps(t)
	\qquad \mbox{for all $t>0$ and } \eps\in (0,1).
  \eas
  As $\sup_{\eps\in (0,1)} \sup_{t>0} \int_{(t-1)_+}^t \heps(s) ds$ is finite thanks to Lemma \ref{lem3}, from Lemma \ref{lem99}
  we directly obtain (\ref{4.1}) with some suitably large $C(p)>0$, because $\frac{p}{p-1}>1$.
\qed
Combining the latter with well-known smoothing estimates for the heat semigroup we can derive uniform boundeness of $\na \veps$.
\begin{lem}\label{lem5}
  Let $n\le 2$. Then there exists $C>0$ such that
  \be{5.1}
	\|\na\veps(\cdot,t)\|_{L^\infty(\Om)} \le C
	\qquad \mbox{for all $t>0$ and } \eps\in (0,1).
  \ee
\end{lem}
\proof
  This follows from Lemma \ref{lem4} when applied to any fixed $p>2$ and combined with (\ref{vinfty}) and standard smoothing 
  estimates for the Neumann heat semingroup on $\Om$ (\cite{win_JDE}).
\qed
This in turn enables us to improve our knowledge on the first solution components:
\begin{lem}\label{lem6}
  When $n\le 2$, one can find $C>0$ with the property that
  \be{6.1}
	\|\ueps(\cdot,t)\|_{L^\infty(\Om)} \le C
	\qquad \mbox{for all $t>0$ and } \eps\in (0,1).
  \ee
\end{lem}
\proof
  On the basis of Lemma \ref{lem4}, Lemma \ref{lem5}, (\ref{phi}) and (\ref{vinfty}), this can be obtained by means of
  a Moser-type iterative argument, as recorded in \cite[Lemma A.1]{taowin_subcrit}.
\qed
As a direct consequence of standard parabolic regularity theory, we next obtain H\"older bounds for $\ueps$ and $\veps$.
\begin{lem}\label{lem7}
  Let $n\le 2$. Then there exist $\theta\in (0,1)$ and $C>0$ such that
  \be{7.1}
	\|\ueps\|_{C^{\theta,\frac{\theta}{2}}(\bom\times [t,t+1])} \le C
	\qquad \mbox{for all $t>0$ and } \eps\in (0,1),
  \ee
  and that
  \be{7.2}
	\|\veps\|_{C^{\theta,\frac{\theta}{2}}(\bom\times [t,t+1])} \le C
	\qquad \mbox{for all $t>0$ and } \eps\in (0,1).
  \ee
\end{lem}
\proof
  In view of the $L^\infty$ estimates provided by Lemma \ref{lem6} and (\ref{vinfty}), both inequalities immediately result
  from a standard result on H\"older regularity in scalar parabolic equations (\cite{porzio_vespri}), again because
  (\ref{phi}) and (\ref{vinfty}) yield $c_1>0, c_2>0$ and $c_3>0$ such that $c_1\le \phi(\veps) \le c_2$ and
  $|\phi'(\veps)|\le c_3$ in $\Om\times (0,\infty)$ for all $\eps\in (0,1)$.
\qed
The derivation of higher order H\"older bounds, locally in time, thereupon becomes straightforward as well.
\begin{lem}\label{lem8}
  If $n\le 2$, then for all $\tau>0$ and any $T>\tau$ there exist $\theta=\theta(\tau,T)\in (0,1)$ and $C(\tau,T)>0$ such that
  \be{8.1}
	\|\ueps\|_{C^{2+\theta,\frac{\theta}{2}}(\bom\times [\tau,T])} \le C(\tau,T)
	\qquad \mbox{for all } \eps\in (0,1)
  \ee
  and
  \be{8.2}
	\|\veps\|_{C^{2+\theta,\frac{\theta}{2}}(\bom\times [\tau,T])} \le C(\tau,T)
	\qquad \mbox{for all } \eps\in (0,1).
  \ee
\end{lem}
\proof
  Based on Lemma \ref{lem7}, the estimate in (\ref{8.2}) follows from parabolic Schauder theory (\cite{LSU}) applied to the
  second equation in (\ref{0eps}).
  Having (\ref{8.2}) at hand, we can thereupon employ the same token again to readily derive (\ref{8.1}), again relying on
  (\ref{phi}) and (\ref{vinfty}).
\qed
Our main result on the existence of global classical solutions to the one- and two-dimensional versions of 
(\ref{0}), as well as their large time behavior, can now be established by using the boundedness properties collected above
with the approximation and stabilization results known from Lemma \ref{lem1}.\abs
\proofc of Theorem \ref{theo9}. \quad
  As a consequence of Lemma \ref{lem7}, Lemma \ref{lem8} and the Arzel\`a-Ascoli theorem, possibly after modification of
  the functions $u$ and $v$ from Lemma \ref{lem1} on a null set we can achieve that in (\ref{1.1}), with $(\eps_j)_{j\in\N}$
  as provided there we have
  \be{9.3}
	\ueps \to u
	\quad \mbox{and } \veps\to v 
	\quad \mbox{in $C^0_{loc}(\bom\times [0,\infty))$ and in $C^{2,1}_{loc}(\bom\times (0,\infty))$ \quad as } 
		\eps=\eps_j\searrow 0,
  \ee
  and that hence $(u,v)$ has the regularity features in (\ref{9.01}) and solves (\ref{0}) classically.
  Since (\ref{9.3}) together with (\ref{7.1}) ensures that
  $(u(\cdot,t))_{t>0}$ is bounded in $C^\theta(\bom)$ with some $\theta\in (0,1)$, 
  and since again due to the boundedness we know that the first of the two continuous embeddings 
  $C^\theta(\bom) \hra L^\infty(\Om)\hra (W^{1,2}(\Om))^\star$ is compact,
  once more drawing on an Ehrling type inequality we infer that writing $\ouz:=\frac{1}{|\Om|} \io u_0$ we see that for each 
  $\eta>0$ there exists $c_1(\eta)>0$ fulfilling
  \bas
	\limsup_{t\to\infty} \|u(\cdot,t)-\ouz\|_{L^\infty(\Om)}
	&\le& \eta \sup_{t>0} \|u(\cdot,t)-\ouz\|_{C^\theta(\bom)}
	+ c_1(\eta) \limsup_{t\to\infty} \|u(\cdot,t)-\ouz\|_{(W^{1,2}(\Om))^\star} \\
	&=& \eta \sup_{t>0} \|u(\cdot,t)-\ouz\|_{C^\theta(\bom)}.
  \eas
  Combined with a similar argument for the second solution component, this shows that also (\ref{9.1}) is valid.
\qed
\mysection{The case $n\ge 3$. Proof of Theorem \ref{theo17}}
\subsection{A space-time $L^\frac{4}{3}$ estimate for $\na\ueps$}
We now turn our attention to the case $n\ge 3$, in which due to lacking favorable embeddings,
any expedient analogue of Lemma \ref{lem4} seems absent. 
To show that nevertheless the pair $(u,v)$ obtained in Lemma \ref{lem1} solves (\ref{0}) actually in a standard
weak sense,		
we use the $L^4$ estimate for $\na\veps$ from Lemma \ref{lem3} in the context of an $L\log L$ testing procedure,
indeed leading to bounds for $\na \ueps$ in some reflexive $L^p$ space.
\begin{lem}\label{lem10}
  Let $n\ge 1$. Then there exists $C\ge 0$ such that
  \be{10.1}
	\io \ueps(\cdot,t) \ln\ueps(\cdot,t) \le C
	\qquad \mbox{for all $t>0$ and } \eps\in (0,1),
  \ee
  and that
  \be{10.2}
	\int_{(t-1)_+}^t \io |\na\ueps|^\frac{4}{3} \le C
	\qquad \mbox{for all $t>0$ and } \eps\in (0,1).
  \ee
\end{lem}
\proof
  Again by (\ref{0eps}) and Young's inequality,
  \bas
	\frac{d}{dt} \io \ueps\ln\ueps
	&=& - \io \phi(\veps) \frac{|\na\ueps|^2}{\ueps}
	- \io \phi'(\veps) \na\ueps\cdot\na\veps \\
	&\le& - c_1 \io \frac{|\na\ueps|^2}{\ueps} 
	+ c_2 \io |\na\ueps|\cdot |\na\veps| \\
	&\le& - \frac{c_1}{2} \io \frac{|\na\ueps|^2}{\ueps}
	+ \frac{c_2^2}{2c_1} \io \ueps |\na\veps|^2 \\
	&\le& - \frac{c_1}{2} \io \frac{|\na\ueps|^2}{\ueps}
	+ \frac{c_2^2}{4c_1} \io \ueps^2
	+ \frac{c_2^2}{4c_1} \io |\na\veps|^4
	\qquad \mbox{for all $t>0$ and } \eps\in (0,1),
  \eas
  where $c_1:=\inf_{\eps\in (0,1)} \inf_{\Om\times (0,\infty)} \phi(\veps)>0$ and  
  $c_2:\sup_{\eps\in (0,1)} \sup_{\Om\times (0,\infty)} |\phi'(\veps)|<\infty$
  by (\ref{phi}) and (\ref{vinfty}).
  As furthermore
  \bas
	- \frac{|\Om|}{e} \le \io \ueps\ln\ueps \le \io \ueps^2
	\qquad \mbox{for all $t>0$ and } \eps\in (0,1)
  \eas
  due to the fact that $-\frac{1}{e} \le \xi\ln\xi\le \xi^2$ for all $\xi>0$, we thus infer that
  \bas	
	\yeps(t):=\io \ueps(\cdot,t)\ln\ueps(\cdot,t) + \frac{|\Om|}{e},
	\qquad t\ge 0, \ \eps\in (0,1),
  \eas
  as well as
  \bas
	& & \hs{-20mm}
	\geps(t):=\frac{c_1}{2} \io \frac{|\na\ueps(\cdot,t)|^2}{\ueps(\cdot,t)}
	\quad \mbox{and} \quad
	\heps(t):=\Big(\frac{c_2^2}{4c_1} +1\Big) \io \ueps^2(\cdot,t) 
	+ \frac{c_2^2}{4c_1} \io \big|\na\veps(\cdot,t)\big|^4 
	+ \frac{|\Om|}{e}, \\[1mm]
	& & \hs{100mm}
	\qquad t>0, \ \eps\in (0,1),
  \eas
  are all nonnegative and satisfy
  \bas
	\yeps'(t) + \yeps(t) + \geps(t) \le \heps(t)
	\qquad \mbox{for all $t>0$ and } \eps\in (0,1).
  \eas
  Since $c_3:=\sup_{\eps\in (0,1)} \sup_{t>0} \int_{(t-1)_+}^t \heps(s) ds$ is finite according to Lemma \ref{lem2} and
  Lemma \ref{lem3}, again in view of Lemma \ref{lem_JFA} this implies that
  \be{10.3}
	\yeps(t) \le \frac{c_3}{1-e^{-1}} + \io u_0^2 + \frac{|\Om|}{e} 
	\qquad \mbox{for all $t>0$ and } \eps\in (0,1),
  \ee
  and that hence
  \be{10.4}
	\int_{(t-1)_+}^t \io \frac{|\na\ueps|^2}{\ueps}
	= \frac{2}{c_1} \int_{(t-1)_+}^t \geps(s) ds
	\le c_4:= \frac{2}{c_1} \cdot \Big( \frac{c_3}{1-e^{-1}} +c_3 + \io u_0^2 + \frac{|\Om|}{e}\Big)
  \ee
  for all $t>0$ and $\eps\in (0,1)$.
  While (\ref{10.3}) directly yields (\ref{10.1}), to derive (\ref{10.2}) it is sufficient to employ Young's inequality in 
  estimating
  \bas
	\int_{(t-1)_+}^t \io |\na\ueps|^\frac{4}{3}
	&=& \int_{(t-1)_+}^t \io \Big( \frac{|\na\ueps|^2}{\ueps}\Big)^\frac{2}{3} \cdot\ueps^\frac{2}{3} \\
	&\le& \int_{(t-1)_+}^t \io \frac{|\na\ueps|^2}{\ueps}
	+ \int_{(t-1)_+}^t \io \ueps^2
	\qquad \mbox{for all $t>0$ and } \eps\in (0,1),
  \eas
  and to combine (\ref{10.4}) with Lemma \ref{lem2}.
\qed
\subsection{Eventual bounds in the case $n=3$}
Of crucial importance for our derivation of the eventual regularity features claimed in Theorem \ref{theo17} will
be the information on asymptotic smallness of $v$ contained in (\ref{1.2}).
Together with an $L^\infty$ approximation property implied by Lemma \ref{lem3} in the three-dimensional case,
this entails a doubly uniform ultimate smallness property of $\veps$ in the following sense.
\begin{lem}\label{lem12}
  Let $n= 3$, and let $(\eps_j)_{j\in\N}$ and $(u,v)$ be as in Lemma \ref{lem1}.
  Then there exists a subsequence $(\eps_{j_k})_{k\in\N}$ of $(\eps_j)_{j\in\N}$ with the property that whenever $\eta>0$,
  one can find $T(\eta)>0$b and $\eps_\star(\eta)\in (0,1)$ such that
  \be{12.1}
	\|\veps(\cdot,t)\|_{L^\infty(\Om)} \le \eta
	\qquad \mbox{for all $t>T(\eta)$ and } \eps\in (\eps_{j_k})_{k\in\N} \subset (0,\eps_\star(\eta)).
  \ee
\end{lem}
\proof
  From Lemma \ref{lem3} it follows that $(\veps)_{\eps\in (0,1)}$ is bounded in $L^4((0,T);W^{1,4}(\Om))$ and that
  $(v_{\eps t})_{\eps\in (0,1)}$ is bounded in $L^2(\Om\times (0,T))$ for all $T>0$, so that since $W^{1,4}(\Om)$ is compactly
  embedded into $L^\infty(\Om)$ due to our assumption that $n\le 3$, an Aubin-Lions lemma (\cite{temam}) along with
  (\ref{1.1}) asserts that $\veps\to v$ in $L^4_{loc}([0,\infty);L^\infty(\Om))$ as $\eps=\eps_j\searrow 0$;
  for an appropriate null set $N_1\subset (0,\infty)$ and some subsequence $(\eps_{j_k})_{k\in\N}$ of $(\eps_j)_{j\in\N}$,
  it thus follows that
  \be{12.2}
	\veps(\cdot,t)\to v(\cdot,t)
	\quad \mbox{in $L^\infty(\Om)$ \quad for all $t\in (0,\infty)\sm N_1$ \qquad as } \eps=\eps_{j_k}\searrow 0.
  \ee
  Now if we fix $\eta>0$ and let $N\subset (0,\infty)$ be as in Lemma \ref{lem1}, in line with (\ref{1.2}) we can pick 
  $T(\eta)\in (0,\infty)\sm (N\cup N_1)$ such that $\|v(\cdot,T(\eta))\|_{L^\infty(\Om)} \le \frac{\eta}{2}$,
  and since then $T(\eta)\in (0,\infty)\sm N_1$, we may rely on (\ref{12.2}) to choose $\eps_\star(\eta)\in (0,1)$
  in such a way that $\|\veps(\cdot,T(\eta))-v(\cdot,T(\eta))\|_{L^\infty(\Om)} \le \frac{\eta}{2}$ for all
  $\eps\in (\eps_{j_k})_{k\in\N}$ such that $\eps<\eps_\star(\eta)$.
  Consequently, 
  \bas
	\big\|\veps(\cdot,T(\eta))\big\|_{L^\infty(\Om)}
	&\le& \big\|\veps(\cdot,T(\eta))-v(\cdot,T(\eta))\big\|_{L^\infty(\Om)}
	+ \big\|v(\cdot,T(\eta))\big\|_{L^\infty(\Om)} \nn\\
	&\le& \frac{\eta}{2}+\frac{\eta}{2}=\eta
	\qquad \mbox{for all } \eps\in (\eps_{j_k})_{k\in\N} \cap (0,\eps_\star(\eta)),
  \eas
  so that (\ref{12.1}) results from the monotonocity property of $0\le t \mapsto \|\veps(\cdot,t)\|_{L^\infty(\Om)}$
  for $\eps\in (0,1)$, as expressed in (\ref{vinfty}).
\qed
By means of an argument dating back to \cite{taowin_consumption},
the key step toward our large time analysis in the three-dimensional case can now be accomplished
through the study of the evolution of the functionals $\io \ueps^p (\del-\veps)^{-\kappa}$ with arbitrary $p>1$
and suitably chosen $\del=\del(p)>0$ and $\kappa=\kappa(p)>0$, 
crucially based on the premise that $\veps$ is sufficiently small.
\begin{lem}\label{lem13}
  Let $n=3$, and let $(\eps_{j_k})_{k\in\N}$ be as provided by Lemma \ref{lem12}.
  Then for each $p>1$ there exist $T(p)>0, \eps_\star(p)\in (0,1)$ and $C(p)>0$ such that
  \be{13.1}
	\io \ueps^p(\cdot,t) \le C(p)
	\qquad \mbox{for all $t>T(p)$ and } \eps\in (\eps_{j_k})_{k\in\N} \subset (0,\eps_\star(p)).
  \ee
\end{lem}
\proof
  Once more relying on (\ref{phi}) and (\ref{vinfty}), we fix $c_1>0, c_2>0$ and $c_3>0$ such that
  \be{13.2}
	c_1 \le \phi(\veps)\le c_2 
	\quad \mbox{and} \quad
	|\phi'(\veps)|\le c_3
	\quad \mbox{in $\Om\times (0,\infty)$ \qquad for all } \eps\in (0,1),
  \ee
  and given $p>1$ we then pick $\kappa=\kappa(p)>0$ such that
  \be{13.4}
	\frac{p(c_2+1)^2\kappa^2}{2(p-1)c_1} \le \frac{\kappa}{2}
  \ee
  and choose $\del=\del(p)>0$ small enough fulfilling both
  \be{13.3}
	p c_3 \del \le \frac{1}{2}
	\qquad \mbox{and} \qquad
	(p-1) c_3\del + \kappa\del \le\kappa.
  \ee
  Thanks to Lemma \ref{lem12}, it is then possible to find $T_1=T_1(p)>0$ and $\eps_\star=\eps_\star(p)\in (0,1)$ such that
  \be{13.5}
	\veps(x,t)\le\frac{\del}{2}
	\qquad \mbox{for all $x\in\Om, t>T_1$ and } \eps\in (\eps_{j_k})_{k\in\N} \cap (0,\eps_\star),
  \ee
  and to make appropriate use of this, we recall (\ref{0eps}) to compute
  \bea{13.6}
	& & \hs{-16mm}
	\frac{d}{dt} \io \ueps^p (\del-\veps)^{-\kappa} \nn\\
	&=& p \io \ueps^{p-1} (\del-\veps)^{-\kappa} \Del \big(\ueps\phi(\veps)\big)
	+ \kappa \io \ueps^p (\del-\veps)^{-\kappa-1} \cdot \Big\{ \Del\veps-\frac{\ueps\veps}{1+\eps\ueps}\Big\} \nn\\
	&=& - p\io \Big\{ (p-1)\ueps^{p-2} (\del-\veps)^{-\kappa} \na\ueps
	+ \kappa\ueps^{p-1} (\del-\veps)^{-\kappa-1} \na\veps\Big\}
		\cdot \Big\{ \phi(\veps)\na\ueps + \ueps\phi'(\veps)\na\veps \Big\} \nn\\
	& & - \kappa \io \Big\{ p\ueps^{p-1} (\del-\veps)^{-\kappa}\na\ueps + (\kappa+1)\ueps^p (\del-\veps)^{-\kappa-2} \na\veps
		\Big\} \cdot\na\veps \nn\\
	& & - \kappa \io \frac{\ueps^{p+1}\veps (\del-\veps)^{-\kappa-1}}{1+\eps\ueps} \nn\\
	&=& -p(p-1) \io \ueps^{p-2} (\del-\veps)^{-\kappa} \phi(\veps) |\na\ueps|^2 \nn\\
	& & - \io \ueps^p \cdot \Big\{ \kappa(\kappa+1) (\del-\veps)^{-\kappa-2} + p\kappa (\del-\veps)^{-\kappa-1} \phi'(\veps)
		\Big\} |\na\veps|^2 \nn\\
	& & - \io \ueps^{p-1} \cdot \Big\{ p\kappa (\del-\veps)^{-\kappa-1} \phi(\veps) 
		+ p(p-1) (\del-\veps)^{-\kappa} \phi'(\veps) + p\kappa (\del-\veps)^{-\kappa} \Big\} \na\ueps\cdot\na\veps \nn\\
	& & - \kappa \io \frac{\ueps^{p+1} \veps (\del-\veps)^{-\kappa-1}}{1+\eps\ueps}
 	\qquad \mbox{for all $t>T_1$ and } \eps\in (\eps_{j_k})_{k\in\N} \cap (0,\eps_\star).
  \eea
  Here by (\ref{13.2}),
  \bea{13.7}
	& & \hs{-20mm}
	p(p-1) \io \ueps^{p-2} (\del-\veps)^{-\kappa} \phi(\veps) |\na\ueps|^2 \nn\\
	&\ge& p(p-1) c_1 \io \ueps^{p-2} (\del-\veps)^{-\kappa} |\na\ueps|^2
	\qquad \mbox{for all $t>T_1$ and } \eps\in (\eps_{j_k})_{k\in\N} \cap (0,\eps_\star),
  \eea
  while due to (\ref{13.2}) and (\ref{13.3}),
  \bas
	& & \hs{-30mm}
	\kappa(\kappa+1) (\del-\veps)^{-\kappa-2} + p\kappa (\del-\veps)^{-\kappa-1} \phi'(\veps) \\
	&\ge& \kappa(\kappa+1) (\del-\veps)^{-\kappa-2} - p\kappa c_3 (\del-\veps)^{-\kappa-1} \\
	&=& \kappa (\del-\veps)^{-\kappa-2} \cdot \big\{ \kappa+1- pc_3\cdot (\del-\veps) \big\} \\
	&\ge& \kappa (\del-\veps)^{-\kappa-2} \cdot \big\{ \kappa+1- pc_3 \del \big\} \\
	&\ge& \frac{\kappa}{2} (\del-\veps)^{-\kappa-2}
	\quad \mbox{in } \Om\times (T_1,\infty)
	\qquad \mbox{for all } \eps\in (\eps_{j_k})_{k\in\N} \cap (0,\eps_\star)
  \eas
  and hence
  \bea{13.8}
	& & \hs{-20mm}
	\io \ueps^p \cdot \Big\{ \kappa(\kappa+1) (\del-\veps)^{-\kappa-2} + p\kappa (\del-\veps)^{-\kappa-1} \phi'(\veps)
		\Big\} |\na\veps|^2 \nn\\
	&\ge& \frac{\kappa}{2} \io \ueps^p (\del-\veps)^{-\kappa-2} |\na\veps|^2
	\qquad \mbox{for all $t>T_1$ and } \eps\in (\eps_{j_k})_{k\in\N} \cap (0,\eps_\star).
  \eea
  Moreover, again by (\ref{13.2}) and (\ref{13.3}),
  \bas
	& & \hs{-20mm}
	\Big| p\kappa (\del-\veps)^{-\kappa-1} \phi(\veps) 
		+ p(p-1) (\del-\veps)^{-\kappa} \phi'(\veps) + p\kappa (\del-\veps)^{-\kappa} \Big| \\
	&\le& p\kappa c_2 (\del-\veps)^{-\kappa-1} + p(p-1) c_3 (\del-\veps)^{-\kappa} + p\kappa (\del-\veps)^{-\kappa} \\
	&=& p (\del-\veps)^{-\kappa-1} \cdot \big\{\kappa c_2 + (p-1) c_3\cdot (\del-\veps) + \kappa\cdot (\del-\veps) \big\} \\
	&\le& p (\del-\veps)^{-\kappa-1} \cdot \big\{ \kappa c_2 + (p-1) c_3 \del + \kappa\del \big\} \\
	&\le& p(c_2+1) \kappa (\del-\veps)^{-\kappa-1}
	\quad \mbox{in } \Om\times (T_1,\infty)
	\qquad \mbox{for all } \eps\in (\eps_{j_k})_{k\in\N} \cap (0,\eps_\star),
  \eas
  so that thanks to Young's inequality,
  for all $t>T_1$ and $\eps\in (\eps_{j_k})_{k\in\N} \cap (0,\eps_\star)$ we have
  \bea{13.9}
	& & \hs{-20mm}
	- \io \ueps^{p-1} \cdot \Big\{ p\kappa (\del-\veps)^{-\kappa-1} \phi(\veps) 
		+ p(p-1) (\del-\veps)^{-\kappa} \phi'(\veps) + p\kappa (\del-\veps)^{-\kappa} \Big\} \na\ueps\cdot\na\veps \nn\\
	&\le& p(c_2+1)\kappa \io \ueps^{p-1} (\del-\veps)^{-\kappa-1} |\na\ueps| \cdot |\na\veps| \nn\\
	&\le& \frac{p(p-1) c_1}{2} \io \ueps^{p-2} (\del-\veps)^{-\kappa} |\na\ueps|^2
	+ \frac{p(c_2+1)^2 \kappa^2}{2(p-1) c_1} \io \ueps^p (\del-\veps)^{-\kappa-2} |\na\veps|^2 \nn\\
	&\le& \frac{p(p-1) c_1}{2} \io \ueps^{p-2} (\del-\veps)^{-\kappa} |\na\ueps|^2
	+ \frac{\kappa}{2} \io \ueps^p (\del-\veps)^{-\kappa-2} |\na\veps|^2
  \eea
  because of (\ref{13.4}).
  From (\ref{13.6})-(\ref{13.9}) we thus obtain, writing 
  \bas
	\yeps(t):=\io \ueps^p(\cdot,t) \big(\del-\veps(\cdot,t)\big)^{-\kappa},
	\qquad t\ge T_1, \eps\in (0,1),
  \eas
  that
  \be{13.10}
	\yeps'(t) + \frac{p(p-1) c_1}{2} \io \ueps^{p-2} (\del-\veps)^{-\kappa} |\na\ueps|^2
	\le 0
	\qquad \mbox{for all $t>T_1$ and } \eps\in (\eps_{j_k})_{k\in\N} \cap (0,\eps_\star),
  \ee
  and here a superlinear absorptive summand can be created by observing that due to the Gagliardo-Nirenberg inequality and
  (\ref{mass}), using (\ref{13.5}) and abbreviating $a\equiv a(p):=\frac{3(p-1)}{3p-1}\in (0,1)$ we can find $c_4=c_4(p)>0$
  and $c_5=c_5(p)>0$ such that
  \bas
	& & \hs{-20mm}
	\bigg\{ \io \ueps^p (\del-\veps)^{-\kappa} \bigg\}^\frac{1}{a} \\
	&\le& \Big(\frac{\del}{2}\Big)^{-\frac{\kappa}{a}} \|\ueps^\frac{p}{2}\|_{L^2(\Om)}^\frac{2}{a} \\
	&\le& c_4 \|\na\ueps^\frac{p}{2}\|_{L^2(\Om)}^2 \|\ueps^\frac{p}{2}\|_{L^\frac{2}{p}(\Om)}^\frac{2(1-a)}{a}
		+ c_4\|\ueps^\frac{p}{2}\|_{L^\frac{2}{p}(\Om)}^\frac{2}{a} \\
	&\le& c_5 \io \ueps^{p-2} |\na\ueps|^2 + c_5 \\
	&\le& c_5\del^\kappa \io \ueps^{p-2} (\del-\veps)^{-\kappa} |\na\ueps|^2
	+ c_5
	\qquad \mbox{for all $t>T_1$ and } \eps\in (\eps_{j_k})_{k\in\N} \cap (0,\eps_\star).
  \eas
  Therefore, (\ref{13.10}) implies that if we let
  $c_6\equiv c_6(p):=\frac{p(p-1) c_1}{2c_5 \del^\kappa}$ and
  $c_7\equiv c_7(p):=\frac{p(p-1) c_1}{2\del^\kappa}$, then
  \bas
	\yeps'(t) + c_6 \yeps^\frac{1}{a}(t) \le c_7
	\qquad \mbox{for all $t>T_1$ and } \eps\in (\eps_{j_k})_{k\in\N} \cap (0,\eps_\star),
  \eas
  so that since $\frac{1}{a}>1$, Lemma \ref{lem122} applies so as to assert that
  \bas
	\yeps(t)
	&\le& \Big(\frac{c_7}{c_6}\Big)^a
	+ \Big(\frac{a}{c_6(1-a)}\Big)^\frac{a}{1-a} \cdot (t-T_1)^{-\frac{a}{1-a}} \\
	&\le& \Big(\frac{c_7}{c_6}\Big)^a
	+ \Big(\frac{a}{c_6(1-a)}\Big)^\frac{a}{1-a}
	\qquad \mbox{for all $t>T_1+1$ and } \eps\in (\eps_{j_k})_{k\in\N} \cap (0,\eps_\star),
  \eas
  and that thus (\ref{13.1}) holds with $T(p):=T_1+1$ upon an evident choice of $C(p)$.
\qed
Again, an application to suitably large $p$ yields pointwise bounds for $\na\veps$.
\begin{lem}\label{lem14}
  Let $n=3$, and let $(\eps_{j_k})_{k\in\N}$ be as given by Lemma \ref{lem12}.
  Then there exist $T>0, \eps_\star\in (0,1)$ and $C>0$ such that
  \be{14.1}
	\|\veps(\cdot,t)\|_{W^{1,\infty}(\Om)} \le C
	\qquad \mbox{for all $t>T$ and } \eps\in (\eps_{j_k})_{k\in\N} \subset (0,\eps_\star).
  \ee
\end{lem}
\proof
  Again based on known regularization features of the Neumann heat semigroup on $\Om$, this can be derived from Lemma \ref{lem13}
  and (\ref{vinfty}) in a straightforward manner.
\qed
After an adequate $\eps$-independent waiting time, both $L^\infty$ and H\"older estimates for $\ueps$ can be obtained
from Lemma \ref{lem13} and Lemma \ref{lem14} upon a suitable cut-off procedure with respect to the time variable.
\begin{lem}\label{lem15}
  Let $n=3$, and let $(\eps_{j_k})_{k\in\N}$ be taken from Lemma \ref{lem12}.
  Then one can find $T>0, \eps_\star\in (0,1), \theta\in (0,1)$ and $C>0$ such that
  \be{15.1}
	\|\ueps\|_{C^{\theta,\frac{\theta}{2}}(\bom\times [t,t+1])} \le C
	\qquad \mbox{for all $t>T$ and } \eps\in (\eps_{j_k})_{k\in\N} \subset (0,\eps_\star).
  \ee
\end{lem}
\proof
  For arbitrary $T>0$ and $\zeta=\zeta^{(T)}\in C^\infty([0,\infty))$ fulfilling $\zeta\equiv 0$ on $[0,T+\frac{1}{2}]$
  and $\zeta\equiv 1$ on $[T+1,\infty)$ as well as $0\le\zeta\le 1$ and $|\zeta'|\le 4$, the functions given by
  \bas
	\weps(x,t):=\weps^{(T)}(x,t):=\zeta(t)\ueps(x,t),
	\qquad (x,t)\in\bom\times [0,\infty), \ \eps\in (0,1),
  \eas
  satisfy $\weps\equiv 0$ on $\bom\times [0,T+\frac{1}{2}]$ and
  \be{15.2}
	w_{\eps t} = \na\cdot \big(D_\eps(x,t)\na\weps\big) - \na\cdot \big(b_\eps(x,t) \weps\big) + \heps(x,t)
  \ee
  for all $\eps\in (0,1)$, with
  $D_\eps(x,t):=\phi(\veps(x,t))$,
  $b_\eps(x,t):=\phi'(\veps(x,t)) \na\veps(x,t)$
  and $\heps(x,t)\equiv \heps^{(T)}(x,t):=\zeta'(t) \ueps(x,t)$,
  $(x,t)\in \Om\times (0,\infty)$, $\eps\in (0,1)$.
  Here we note that $\inf_{\eps\in (0,1)} \inf_{\Om\times (0,\infty)} D_\eps>0$ by (\ref{phi}) and (\ref{vinfty}), and that
  Lemma \ref{lem13}, Lemma \ref{lem14}, (\ref{phi}) and (\ref{vinfty}) assert that for each $p>1$ we can find 
  $T_1=T_1(p)>0$ and $\eps_\star=\eps_\star(p)\in (0,1)$ such that
  $(\weps)_{\eps\in (\eps_{j_k})_{k\in\N} \cap (0,\eps_\star)}$,
  $(|b_\eps|)_{\eps\in (\eps_{j_k})_{k\in\N} \cap (0,\eps_\star)}$ and
  $(\heps)_{\eps\in (\eps_{j_k})_{k\in\N} \cap (0,\eps_\star)}$ 
  are bounded in $L^\infty((T_1,\infty);L^p(\Om))$.
  Applying this to some suitably large $p$ in the course of a Moser-type iteration in (\ref{15.2}) 
  (\cite[Lemma A.1]{taowin_subcrit}), we hence infer the existence of $T_2>0, \eps_\star\in (0,1)$ and $c_1>0$ such that
  \bas
	\|\weps^{(T_2)}(\cdot,t)\|_{L^\infty(\Om)} \le c_1
	\qquad \mbox{for all $t>T_2$ and } \eps\in (\eps_{j_k})_{k\in\N} \cap (0,\eps_\star),
  \eas
  and that thus, by definition of $(\weps^{(T_2)})_{\eps\in (0,1)}$, 
  \bas
	\|\ueps(\cdot,t)\|_{L^\infty(\Om)} \le c_1
	\qquad \mbox{for all $t>T_2+1$ and } \eps\in (\eps_{j_k})_{k\in\N} \cap (0,\eps_\star).
  \eas
  Based on this, and again on Lemma \ref{lem14} and (\ref{phi}), we may draw on temporally localized parabolic H\"older estimates
  (\cite{porzio_vespri}) to see that, indeed, (\ref{15.1}) holds with $\eps_\star$ as above, $T:=T_2+1$ and some suitably
  large $C>0$.
\qed
Higher order regularity features can thereupon directly be inferred from standard parabolic theory.
\begin{lem}\label{lem16}
  Let $n=3$, and let $(\eps_{j_k})_{k\in\N}$ be as in Lemma \ref{lem12}.
  Then there exist $T>0, \eps_\star\in (0,1), \theta\in (0,1)$ and $C>0$ such that
  \be{16.1}
	\|\ueps\|_{C^{2+\theta,1+\frac{\theta}{2}}(\bom\times [t,t+1])} \le C
	\qquad \mbox{for all $t>T$ and } \eps\in (\eps_{j_k})_{k\in\N} \subset (0,\eps_\star)
  \ee
  as well as
  \be{16.2}
	\|\veps\|_{C^{2+\theta,1+\frac{\theta}{2}}(\bom\times [t,t+1])} \le C
	\qquad \mbox{for all $t>T$ and } \eps\in (\eps_{j_k})_{k\in\N} \subset (0,\eps_\star).
  \ee
\end{lem}
\proof  
  By resorting to standard parabolic Schauder theory again (\cite{LSU}), 
  we readily obtain first (\ref{16.2}) and then (\ref{16.1})
  using Lemma \ref{lem15}, Lemma \ref{lem14} and (\ref{phi}).
\qed
We are now in a position to verify our main result on eventual smoothness and asymptotic behavior in the case $n=3$. \abs
\proofc of Theorem \ref{theo17}. \quad
  Taking $(u,v)$ and $(\eps_{j_k})_{k\in\N}$ as provided by Lemma \ref{lem1} and Lemma \ref{lem12},
  from Lemma \ref{lem2}, Lemma \ref{lem10} and Lemma \ref{lem3} we conclude by means of an Aubin-Lions lemma that
  \begin{eqnarray}
	& & \ueps\wto u
	\qquad \mbox{in } L^2_{loc}(\bom\times [0,\infty)),
	\label{17.3} \\
	& & \ueps\wto u
	\qquad \mbox{in } L^\frac{4}{3}_{loc}([0,\infty);W^{1,\frac{4}{3}}(\Om)),
	\label{17.4} \\
	& & \veps\wsto v
	\qquad \mbox{in } L^\infty(\Om\times (0,\infty),
	\label{17.5} \\
	& & \veps\wto v
	\qquad \mbox{in } L^2_{loc}([0,\infty);W^{2,2}(\Om))
	\qquad \mbox{and} 
	\label{17.6} \\
	& & \veps\to v
	\qquad \mbox{in } L^2_{loc}([0,\infty);W^{1,2}(\Om))
	\quad \mbox{and a.e.~in } \Om\times (0,\infty)
	\label{17.7}
  \end{eqnarray}
  as $\eps=\eps_{j_k}\searrow 0$.
  Besides implying (\ref{reg2}), together with a result on stability of weak $L^1$ convergence with respect to certain
  small nonlinear perturbations of the identity (\cite[Lemma 5.1]{liwin2}) which ensures that
  \bas
	\frac{\ueps\veps}{1+\eps\ueps} \wto uv
	\quad \mbox{in } L^1_{loc}(\bom\times [0,\infty))
	\qquad \mbox{as } \eps=\eps_{j_k}\searrow 0
  \eas
  by (\ref{17.3}), (\ref{17.7}) and (\ref{vinfty}), 
  these approximation properties (\ref{17.3})-(\ref{17.7}) enable us to derive (\ref{wu}) and (\ref{wv}) for each
  $\vp\in C_0^\infty(\bom\times [0,\infty))$ from the corresponding weak formulations associated with (\ref{0eps}) in a
  straightforward fashion.\abs
  For adequately large $T>0$,
  the additional smoothness features in (\ref{17.1}) are consequences of Lemma \ref{lem16} and the Arzel\`a-Ascoli theorem,
  whereas (\ref{9.1}) results from (\ref{1.1}) by interpolation, because if $T$ is sufficiently large, then 
  $(u(\cdot,t))_{t>T}$ and $(v(\cdot,t))_{t>T}$ are both bounded in $C^1(\bom)$ due to Lemma \ref{lem16}.
\qed
\mysection{Appendix}
In this appendix we briefly collect three statements on quantitative bounds for linearly and superlinearly 
damped first order ordinary differential inequalities.\abs
A proof of the following can be found in \cite[Lemma 3.4]{win_JFA}.
\begin{lem}\label{lem_JFA}
  Let $t_0\in\R, T>t_0$, $a>0$ and $b>0$, and suppose that $y\in C^0([t_0,T)) \cap C^1((t_0,T))$ and
  $h\in L^1_{loc}(\R)$ are such that $y\ge 0$ on $(t_0,T)$ and $h\ge 0$ a.e.~on $\R$, that
  \be{J1}
	\int_{(t-1)_+}^t h(s) ds \le b
	\qquad \mbox{for all } t\in (t_0,T),
  \ee
  and that
  \bas
	y'(t) + ay(t) \le h(t)
	\qquad \mbox{for all } t\in (t_0,T).
  \eas
  Then
  \bas
	y(t) \le y(t_0) + \frac{b}{1-e^{-a}}
	\qquad \mbox{for all } t\in [t_0,T).
  \eas
\end{lem}
The following elementary inequality has been recorded in \cite[Lemma 2.2]{win_ct_nasto_sig_dep_mot_log}.
\begin{lem}\label{lem99}
  Let $\lam>1$, and suppose that with some $T>0$, $a>0$ and $b>0$,
  the functions $y\in C^0([0,T)) \cap C^1((0,T))$
  and $h\in L^1_{loc}(\R)$ satisfy $y(t)>0$ for all $t\in [0,T)$ and $h(t)\ge 0$ for a.e.~$t\in \R$, and are such that  
  \bas
	\int_{(t-1)_+}^t h(s) ds \le b
	\qquad \mbox{for all } t\in (t_0,T)
  \eas
  as well as
  \bas
	y'(t)+ ay^\lam(t) \le h(t)y(t)
	\qquad \mbox{for all $t\in (0,T)$.}
  \eas
  Then
  \bas
	y(t) \le \max \Big\{ y(0) e^b \, , \, \big\{ a(\lam-1)\big\}^{-\frac{1}{\lam-1}} e^b \Big\}
	\qquad \mbox{for all } t\in [0,T).
  \eas
\end{lem}
Let us finally provide a brief proof for a variant of Lemma \ref{lem99}, here focusing on an estimate
independent of the behavior near the initial instant.
\begin{lem}\label{lem122}
  Let $\lam>1, a>0, b>0, t_0\in\R$ and $T>t_0$, and suppose that $y\in C^1((t_0,T)) \cap L^\infty((t_0,T))$ be nonnegative
  and such that 
  \be{122.1}
	y'(t) + ay^\lam(t) \le b
	\qquad \mbox{for all } t\in (t_0,T).
  \ee
  Then
  \be{122.2}
	y(t)\le \Big(\frac{b}{a}\Big)^\frac{1}{\lam}
	+ \Big(\frac{1}{a(\lam-1)} \Big)^\frac{1}{\lam-1} \cdot (t-t_0)^{-\frac{1}{\lam-1}}
	\qquad \mbox{for all } t\in (t_0,T).
  \ee
\end{lem}
\proof
  Since
  \bas
	\oy(t):=c_1 + c_2\cdot (t-t_0)^{-\frac{1}{\lam-1}},
	\qquad t>t_0,
  \eas
  with $c_1:=(\frac{b}{a})^\frac{1}{\lam}$ and $c_2:=(\frac{1}{a(\lam-1)})^\frac{1}{\lam-1}$, satisfies
  \bas
	\oy'(t) + a\oy^\lam(t) - b
	&=& - \frac{c_2}{\lam-1} (t-t_0)^{-\frac{\lam}{\lam-1}} 
		+ a\cdot \Big\{ c_1 + c_2\cdot (t-t_0)^{-\frac{1}{\lam-1}} \Big\}^\lam - b \\
	&\ge& - \frac{c_2}{\lam-1} (t-t_0)^{-\frac{\lam}{\lam-1}} 
		+ a c_1^\lam + a c_2^\lam (t-t_0)^{-\frac{\lam}{\lam-1}} -b
	\qquad \mbox{for all } t>t_0
  \eas
  due to the fact that $(\xi+\eta)^\lam \ge \xi^\lam + \eta^\lam$ for all $\xi\ge 0$ and $\eta\ge 0$, and since here
  $a c_1^\lam - b=0$ and $a c_2^\lam - \frac{c_2}{\lam-1}=0$, the inequality in (\ref{122.2}) can readily be derived by means 
  of a comparison argument applied to (\ref{122.1}), because $\oy(t)-y(t)\to +\infty$ as $t\searrow t_0$.
\qed

\bigskip

{\bf Acknowledgement.} \
  The first author was funded by the 
  China Scholarship Council (No. 202006630070). 
  The second author acknowledges support of the {\em Deutsche Forschungsgemeinschaft} (Project No.~462888149).
\small


\begin{thebibliography}{99}
%
\bibitem{ahn_yoon}
  \sc Ahn, J., Yoon, C.: 
  \it Global well-posedness and stability of constant equilibria in parabolic-elliptic chemotaxis systems without gradient sensing.
  \rm Nonlinearity {\bf 32}, 1327-1351 (2019)
\bibitem{burger}
  \sc Burger, M., Lauren\c{c}ot, Ph., Trescases, A.:
  \it Delayed blow-up for chemotaxis models with local sensing.
  \rm J.~London Math.~Soc. {\bf 103}, 1596-1617 (2021)  
\bibitem{desvillettes}
  \sc Desvillettes, L., Kim, Y.-J., Trescases, A., Yoon, C.:
  \it A logarithmic chemotaxis model featuring global existence and aggregation.
  \rm Nonlin.~Anal.~Real World Appl. {\bf 50}, 562-582 (2019)   
\bibitem{DLTW}
  \sc Desvillettes, L., Trescases, A., Lauren\c{c}ot, Ph., Winkler, M.:
  \it Weak solutions to triangular cross diffusion systems modeling chemotaxis with local sensing. 
  \rm In preparation
\bibitem{fu}
  \sc Fu, X., Tang, L.H., Liu, C., Huang, J.D., Hwa, T., Lenz, P.:
  \it Stripe formation in bacterial systems with density-suppresses motility.
  \rm Phys. Rev. Lett. {\bf 108}, 198102 (2012)
\bibitem{fujie_jiang_JDE2020}
  \sc Fujie, K., Jiang, J.: \it Global existence for a kinetic model of pattern formation with density-suppressed motilities
  \rm J.~Differential Eq. {\bf 269}, 5338-5378 (2020)   
\bibitem{fujie_jiang_CVPDE}
  \sc Fujie, K., Jiang, J.: 
  \it Comparison methods for a Keller-Segel-type model of pattern formations with density-suppressed motilities.
  \rm Calc. Var. Partial Differential Equations {\bf 60}, 92 (2021)
\bibitem{fujie_jiang_ACAP2021}
  \sc Fujie, K., Jiang, J.: 
  \it Boundedness of classical solutions to a degenerate Keller-Segel type model with signal-dependent motilities.
  \rm Acta Appl.~Math. {\bf 176}, 3 (2021)
\bibitem{fujie_senba}
  \sc Fujie, K., Senba, T.: \it Global existence and infinite time blow-up of classical solutions 
  to chemotaxis systems of local sensing in higher dimensions.
  \rm Preprint. {\tt arXiv:2102.12080}
\bibitem{henry}
  \sc Henry, D.:
  \it Geometric Theory of Semilinear Parabolic Equations.
  \rm Lecture Notes in Mathematics, Vol. {\bf 840}, (Springer- Verlag, 1981)
\bibitem{jiang_arxiv}
  \sc Jiang, J.: \it Boundedness and exponential stabilization in a parabolic-elliptic Keller-Segel model with 
  signal-dependent motilities for local sensing chemotaxis.
  \rm Preprint. {\tt arXiv:2009.07038}
\bibitem{jiang_laurencot}
  \sc Jiang, J., Lauren\c{c}ot, Ph.: \it Global existence and uniform boundedness in a chemotaxis model 
  with signal-dependent motility.
  \rm J. Differential Equations {\bf 299}, 513-541 (2021)
\bibitem{jin_wang}
  \sc Jin, H.-Y., Wang, Z.-A.:
  \it Critical mass on the Keller-Segel system with signal-dependent motility.
  \rm Proc. Amer. Math. Soc. {\bf 148}, 4855-4873 (2020)
\bibitem{jin_kim_wang}
  \sc Jin, H.-Y., Kim, Y.-J., Wang, Z.-A.:
  \it Boundedness, stabilization, and pattern formation driven by density-suppressed motility.
  \rm SIAM J.~Appl.~Math. {\bf 78}, 1632-1657 (2018)  
\bibitem{LSU}
  \sc Ladyzenskaja, O. A., Solonnikov, V. A., Ural'ceva, N. N.:
  \it Linear and Quasi-Linear Equations of Parabolic Type.
  \rm Amer. Math. Soc. Transl., Vol. {\bf 23}, Providence, RI, 1968
\bibitem{KS1}
  \sc Keller, E.F., Segel, L.A.:
  \it Initiation of slime mold aggregation viewed as an instability.
  \rm J. Theor. Biol. {\bf 26}, 399-415 (1970) 
\bibitem{KS3}
  \sc Keller, E.F., Segel, L.A.: \it Model for chemotaxis.
  \rm J. Theoret. Biol. {\bf 30}, 225-234 (1971)
\bibitem{li_zhao_ZAMP}
  \sc Li, D., Zhao, J.:
  \it Global boundedness and large time behavior of solutions to a chemotaxis-consumption system with signal-dependent motility.
  \rm Z.~Angew.~Math.~Physik {\bf 72}, 57 (2021) 
\bibitem{liwin2}
  \sc Li, G., Winkler, M.: \it Relaxation in a Keller-Segel-consumption system involving signal-dependent motilities.
  \rm Preprint
\bibitem{liu}
  \sc Liu, C., et al.:
  \it Sequential establishment of stripe patterns in an expanding cell population.
  \rm Science {\bf 334}, 238 (2011)
\bibitem{liu_xu}
  \sc Liu, Z., Xu, J.: 
  \it Large time behavior of solutions for density-suppressed motility system in higher dimensions.
  \rm J.~Math.~Anal.~Appl. {\bf 475}, 1596-1613 (2019) 
\bibitem{wenbin_lv_ZAMP}
  \sc Lv, W., Wang, Q.: 
  \it Global existence for a class of chemotaxis systems with signal-dependent motility, indirect signal production 
  and generalized logistic source.
  \rm Z.~Angew.~Math.~Physik {\bf 71}, 53 (2020)
\bibitem{wenbin_lv_EECT}
  \sc Lv, W., Wang, Q.: 
  \it Global existence for a class of Keller-Segel model with signal-dependent motility and general logistic term.
  \rm Evol. Equ. Control Theory {\bf 10}, 25-36 (2021)
\bibitem{wenbin_lv_PROCA}
  \sc Lv, W., Wang, Q.: 
  \it A $n$-dimensional chemotaxis system with signal-dependent motility and generalized logistic source: Global existence 
  and asymptotic stabilization.
  \rm Proc. Roy. Soc. Edinburgh Sect. A {\bf 151}, 821-841 (2021)
\bibitem{porzio_vespri}
   \sc Porzio, M.M., Vespri, V.:
   \it Holder estimates for local solutions of some doubly nonlinear degenerate parabolic equations.
   \rm J.~Differential Equations {\bf 103} (1), 146-178 (1993) 
\bibitem{taowin_consumption}
  \sc Tao, Y., Winkler, M.: 
  \it Eventual smoothness and stabilization of large-data solutions 
  in a three-dimensional chemotaxis system with consumption of chemoattractant.
  \rm J.~Differential Equations {\bf 252}, 2520-2543 (2012)  
\bibitem{taowin_subcrit}
  \sc Tao, Y., Winkler, M.: \it Boundedness in a quasilinear parabolic-parabolic Keller-Segel system 
  with subcritical sensitivity.
  \rm J.~Differential Equations {\bf 252}, 692-715 (2012)
\bibitem{taowin_M3AS}
  \sc Tao, Y., Winkler, M.: 
  \it Effects of signal-dependent motilities in a Keller-Segel-type reaction-diffusion system 
  \rm Math.~Mod.~Meth.~Appl.~Sci. {\bf 27}, 1645-1683 (2017)
\bibitem{temam}
  \sc Temam, R.: \it Navier-Stokes Equations. Theory and Numerical Analysis.
  \rm Stud. Math. Appl., Vol. 2, North-Holland, Amsterdam, 1977
\bibitem{wang_wang}
  \sc Wang, J., Wang, M.: 
  \it Boundedness in the higher-dimensional Keller-Segel model with signal-dependent motility and logistic growth.
  \rm J.~Math.~Phys. {\bf 60}, 011507 (2019)
\bibitem{win_JDE}
  \sc Winkler, M.: \it Aggregation vs. global diffusive behavior in the higher-dimensional Keller-Segel model.
  \rm J.~Differential Equations {\bf 248}, 2889-2905 (2010)
\bibitem{win_JFA}
  \sc Winkler, M: \it A three-dimensional Keller-Segel-Navier Stokes system with logistic source: 
  Global weak solutions and asymptotic stabilization.
  \rm J. Funct. Anal. {\bf 276(1)}, 1339–1401 (2019)
\bibitem{win_sig_dep_mot_cons_2}
  \sc Winkler, M.: 
  \it A quantitative strong parabolic maximum principle and application to a taxis-type migration-consumption model
  involving signal-dependent degenerate diffusion.
  \rm Preprint  
\bibitem{win_sig_dep_mot_cons_general}
  \sc Winkler, M.: 
  \it Application of the Moser-Trudinger inequality in the construction of global solutions to a strongly degenerate migration model.
  \rm Preprint    
\bibitem{win_ct_nasto_sig_dep_mot_log}
  \sc Winkler, M.: 
  \it $L^p$ bounds in the two-dimensional Navier-Stokes system and application to blow-up suppression
  in weakly damped chemotaxis-fluid systems.
  \rm Preprint
\bibitem{win_sig_dep_mot_cons_largetime}
  \sc Winkler, M.: \it Stabilization despite pervasive strong cross-degeneracies in a nonlinear diffusion model for 
  migration-consumption interaction.
  \rm Preprint
\bibitem{yifu_wang}
  \sc Xu, C., Wang, Y.: \it Asymptotic behavior of a quasilinear Keller-Segel system with signal-suppressed motility.
  \rm Calc.~Var. Partial Differential Equations {\bf 60}, 183 (2021)
%
\end{thebibliography}
\end{document}